\documentclass[leqno]{siamltex}
\usepackage{amsmath}
\usepackage{graphicx}
\usepackage{mathrsfs}

\usepackage{amsthm}
\usepackage{float}
\usepackage{amsfonts,amssymb,,graphicx}
\usepackage{dsfont}
\usepackage{pifont}
\usepackage{multirow}
\usepackage{float}
\usepackage{diagbox}
\usepackage{multirow}
\usepackage{multicol}
\usepackage{arydshln}
\usepackage{bbm}
\usepackage[ruled,linesnumbered]{algorithm2e}
\usepackage{booktabs}
\usepackage{hyperref}
\hypersetup{hypertex=true,
	colorlinks=true,
	linkcolor=magenta,
	anchorcolor=blue,
	citecolor=blue}
\usepackage{cite}
\usepackage{tikz}
\usepackage{wrapfig} 
\usepackage{hyperref}
\usepackage{multirow}
\usepackage{threeparttable}
\usepackage{mathtools}
\usepackage{mfirstuc}

\usepackage{amsthm,amsmath,amssymb}
 
\usepackage{mathrsfs}

\usepackage{subfigure}
\usepackage{color}

\usepackage{caption}

\title{Adaptive trajectories sampling  for solving PDEs with deep learning methods
}
\author{Xingyu Chen\footnotemark[1]
	\and Jianhuan Cen\footnotemark[1]
        \and Qingsong Zou\footnotemark[2]
	}
\date{}
\begin{document}
\maketitle
\renewcommand{\thefootnote}{\fnsymbol{footnote}}
\footnotetext[1]{School of Computer Science and Engineering, Sun Yat-sen University, Guangzhou, 510006, China.}
\footnotetext[2]{Corresponding author. School of Computer Science and Engineering, and Guangdong Province Key Laboratory of Computational Science, Sun Yat-sen University, Guangzhou 510006, China. Email: mcszqs@mail.sysu.edu.cn.}

\captionsetup[figure]{labelfont={bf},labelformat={default},labelsep=period,name={Figure}}
\captionsetup[table]{labelfont={bf},labelformat={default},labelsep=period,name={Table}}

\begin{abstract}
In this paper, we propose a new adaptive technique, named {\it adaptive trajectories sampling} (ATS), which is used to select training points for the numerical solution of partial differential equations (PDEs) with deep learning methods. The key feature of the ATS is that all training points are adaptively selected from trajectories that are generated according to a  PDE-related stochastic process.
We incorporate the ATS into three known deep learning solvers for PDEs, namely the adaptive derivative-free-loss method (ATS-DFLM), the adaptive physics-informed neural network method (ATS-PINN), and the adaptive temporal-difference method for forward-backward stochastic differential equations (ATS-FBSTD).
Our numerical experiments demonstrate that the ATS  remarkably improves the computational accuracy and efficiency of the original deep learning solvers for the PDEs. In particular, for some specific high-dimensional PDEs, the ATS can even improve the accuracy of the PINN by two orders of magnitude.
\end{abstract}
\begin{keywords}
   Deep learning, PDEs, Adaptive Sampling, Stochastic process
\end{keywords}

\pagestyle{myheadings}\
\section{Introduction}\
Due to many advantages, such as being free from the curse of dimensionality, no mesh partitioning, and being easy to deal with nonlinear problems, the deep learning method becomes  more and more popular in  solving PDEs, see e.g.\cite{DGM, DRM, PINNs, WAN, PIAN, PFNN, optimal-PINN, Derivative, BSDE-2017, BSDE-2018, FBSNN, DBDP, FBSTD}. 
However, comparing to traditional methods (e.g. the finite element method), deep learning methods also suffer from some weaknesses. For example: there is no theory to guarantee its stability and convergence, and for low-dimensional problems, the deep method often has lower accuracy and higher computational cost. 

In the literature, many techniques have been proposed  to improve the performance of the deep learning method\cite{ML-PINN, PINN fail, self-adaptive loss, hard constraints PINN, gPINN, parallel, allen-cahn, soft attention, PPINN, XPINN, net-PINN1, net-PINN2, DeepXDE, EI-PINN, FI-PINN, ADN, DAS-PINN, Selectnet, Rang, ASS-PINN, RAD, ADLGM, RFM}, one of which is to adaptively sample training points. 
For instances, in 2019, a residual-based adaptive refinement (RAR) method was proposed in \cite{DeepXDE} to improve the efficiency of the well-known physics-informed neural network (PINN) method. The RAR improves the distribution of training points during the training process by sampling more points in the locations where the PDE residual is large. 
In 2021, the paper \cite{EI-PINN} generates the set of training points using the importance sampling based on the probability density function  proportional to the residuals of PDEs.
More recently in 2022, the failure-informed adaptive sampling (FI-PINN) method was proposed in \cite{FI-PINN}. By the FI-PINN, the training points are sampled from  so-called  {\it failure regions} where the residuals of the PDEs are larger than a given tolerance. Also in 2022, the paper \cite{ADN} updates the set of training points by adding sampling points to the regions where the residual is relatively large.
In a very recent paper \cite{DAS-PINN}, a so-called DAS-PINN method uses  a residual-based generative model to generate new  training points for further training. On other adaptive sampling techniques, we refer to \cite{Selectnet, Rang, ASS-PINN, RAD, ADLGM, RFM}.

Most of the aforementioned adaptive sampling techniques select training points of the current training step based on the previous step's residual. Here the function of the residual is to measure the (local) error between the exact solution and the trained neural network (NN) solution of the previous step. Namely, the residual is used as a so-called {\it error indicator}, a computable quantity which equivalents more or less to the incalculable error between the exact $u$ and the NN solution $u_\theta$. Since it involves the computation of derivatives, the computation of the residual might be very expensive, especially in the case of high-dimensional PDEs. Moreover, since no theory guarantees the equivalence between the residual indicator and the true error, the residual-type adaptive sampling might not be very efficient either.

In this paper, we construct a novel adaptive sampling technique completely different from the above residual type method.
Firstly, all our training points are obtained from trajectories that are generated according to a PDE-related stochastic process. The advantage of sampling in this way is that the training points of the two successive steps are closely related, and we can fully use the information from the training points of the previous training step for the current training step.
Secondly, the error indicator in this paper is designed in a way distinct from the residual. In fact, to calculate the error $|u-u_\theta|$ at a certain sampling point $x_0$, we use the so-called {\it empirical value} of $x_0$ to approximate the incalculable $u(x_0)$. 
Here, the empirical value at $x_0$ is often calculated by using the values of $u_\theta$ at some neighbor sampling points which  are generated by using the aforementioned  PDE-related stochastic process, starting from $x_0$. Typically, we can choose the empirical value at $x_0$ to be the average of $u_\theta$ at its neighbor sampling points plus some {\it rewards}, which is a PDE-related function to record the information by moving from $x_0$ to its neighbor points. One advantage of this kind of error indicator is its computation does not involve derivatives and thus it greatly reduces the computational cost of the adaptive sampling. 

We call the above adaptive sampling as the {\it adaptive trajectories sampling} (ATS) technique. Note that in the ATS, both the generation and selection of training points are concerned with a PDE-related stochastic process. With the stochastic process, the information between two successive steps is combined together to guide the sampling of training points. Fundamentally, the basic idea of the ATS is originated from the {\it temporal-difference} idea in the reinforcement learning\cite{RL}.


We would like to mention that the ATS is independent of the original deep learning solvers for the PDEs. In other words, the ATS can be applied to improve the performance of any deep solvers for the PDEs. In fact, we incorporate the ATS into three well-known deep learning solvers for PDEs to construct: the adaptive derivative-free loss method (ATS-DFLM), the adaptive PINN method (ATS-PINN), and the adaptive temporal-difference method for forward-backward stochastic differential equations (ATS-FBSTD).
Our numerical results show that compared with original deep learning methods, the ATS markedly improves the computational accuracy with a small additional computational cost. In particular, for some high-dimensional cases, the relative errors by the ATS-PINN can achieve the order of $O(10^{-4})$, in comparison to the $O(10^{-2})$ by the PINN itself.

The rest of the paper is organized as the following.
In our core Section 2,  we introduce the basic idea of the  ATS technique. In Section 3, we apply the ATS to the DFLM method to design the adaptive deep learning method ATS-DFLM. In Sections 4 and 5 we apply the ATS to the deep learning methods PINN and FBSTD to construct the adaptive deep learning methods ATS-PINN and ATS-FBSTD, respectively. In Section 6, we present some numerical experiments to illustrate how the ATS improves the performance of deep learning solvers for PDEs.

\section{The adaptive trajectories sampling(ATS)}\ 

We consider the following partial differential equations(PDEs)
\begin{equation}\label{PDE}
\begin{cases}
    \mathcal{Q}(u):=\mathcal{L}(u)-f=0, \ & {\rm in}\ \Omega,\\
    \mathcal{B}(u)=g, \ &{\rm on \ } \partial \Omega,
\end{cases}
\end{equation}
where $u$ is an unknown function defined in the bounded domain $\Omega\subset {\mathbb R}^d$,
$f$ and $g$ are given functions on $\Omega$ and $\partial\Omega$. Respectively, $\mathcal{L}$ is an interior differential operator,  $\mathcal{B}$ is  a boundary  operator such as $\mathcal{B} = I$ or $\mathcal{B} = \frac{\partial u}{\partial n}$.

As a standard deep learning  method, we use an {\it neural network}(NN) function to simulate the exact solution $u$ of \eqref{PDE}. For instance, we may simulate $u$ with a residual neural network (ResNet, \cite{Res-net}) whose architecture consists of an input layer, an output layer, and several residual connection blocks  which  take the form 
\[
\begin{array}{c}
\mathbf{y}_{0}:=\delta\left(\mathbf{W}_{0} \mathbf{x}+\mathbf{b}_{0}\right), \\

\mathbf{y}_{i+1}:=\mathbf{B}(\mathbf{y}_i), i=0, \ldots, N-1,\\

u_\theta(\mathbf{x}):=\mathbf{W}_{N} \mathbf{y}_{N}+\mathbf{b}_{N},\\
\end{array}
\]
in which $\mathbf{B}$ is a residual block defined as
\[
\begin{array}{c}
\mathbf{y}_{i,0}:=\mathbf{y}_i,\\

\mathbf{y}_{i,j+1}=\delta\left(\mathbf{W}_{i,j}\mathbf{y}_{i,j}+\mathbf{b}_{i,j}\right), j=0,\dots,M-1,\\

\mathbf{B}(\mathbf{y}_i):=\mathbf{y}_i + \mathbf{y}_{i,M}.
\end{array}
\]
Here, $x\in \mathbb{R}^n$ is the input variable, $u_\theta(x)$ is the corresponding output, and $\delta:\mathbb{R}\rightarrow \mathbb{R}$ is the activation function applied element-wisely to a vectorial function. Moreover, we denote the  parameters by $\boldsymbol{\theta}$, i.e.,
\[
    \boldsymbol{\theta}=\left\{\mathbf{W}_{0}, \mathbf{b}_{0}, \mathbf{W}_{i,j}, \mathbf{b}_{i,j}, \mathbf{W}_{N}, \mathbf{b}_{N} \mid 0 \leq i \leq N-1, 0 \leq j \leq M-1\right\}.
\]
$N$ represents the number of residual blocks of the ResNet and $M$ represents the depth of the residual block. To well simulate $u$ with $u_\theta$, we should train the parameters $\boldsymbol{\theta}$ using a suitable {\it loss function}.

Usually, a deep PDE solver trains the above NN solution by randomly sampling training points. When the solution of a PDE is not uniformly  regular, random sampling may lead to inadequate training  in these singular or highly oscillatory local regions. In order to improve the accuracy of the NN solution, it is necessary to sample more training points in those singular or highly oscillatory regions. However, since the exact solution $u$ is unknown in advance, the position of its singularity/high-frequency oscillation is also unknown. Therefore, usually one often uses the information of the NN solution obtained in the previous training step to approximately locate the singular or high oscillation position of the solution. This method, which  samples more training points in the singular region located according to the previous training solution, is commonly known as the adaptive sampling technique.

In this paper, we introduce a novel adaptive sampling technique: the so-called {\it adaptive trajectories sampling}(ATS) technique. One important feature of this sampling technique is to keep the number of sampling points for each training unchanged so that we can obtain  a fixed number of sampling trajectories by connecting the sampling points of the two successive training steps with line segments. Another feature of this technique is that the sampling points of the current training step are generated adaptively according to the previous training step's sampling points and the trained NN solution. Below we explain the details of our adaptive trajectories sampling.

\subsection{Generate admissible sampling points with a stochastic process}
In this subsection, we generate admissible sampling points for the current training step based on the sampling points of the previous training step.

Suppose we will train the NN function $N$ steps which are denoted by $n=0,1,\ldots, N-1$, we denote their corresponding set of training points by ${\cal S}_n$. By our previous explanation, we fix the cardinality of ${\cal S}_n$ to be a prescribed number $I>0$. In the following, we explain how to generate ${\cal S}_n$ for $n=0,1,\ldots, N-1$.

First, we let ${\cal S}_0$ be the set of points sampled from $\Omega$ according to a certain distribution(e.g. the {\it Uniform} or the {\it Guassian}). For $n\ge 1$, we do not 
generate ${\cal S}_n$ directly, but first generate  ${\cal S}_n'$, a set of admissible sampling points which has a larger cardinality than ${\cal S}_n$,
 based on  $\mathcal{S}_{n-1}=\{{\bf x}_i|i=1, \ldots, I\}$, the set of training points in the $(n-1)$th step.  Below we explain how to generate $S'={\cal S}_n'$ from  $S=\mathcal{S}_{n-1}$.

Let $\{X_t:t\in[0,\infty)\}\in \mathbb{R}^d$ be a $d$-dimensional stochastic process, 
then for any $0<t_1<t_2$, we have
\begin{equation}\label{stochastic process}
    X_{t_2}=X_{t_1}+\int_{t_1}^{t_2} d X_t.
\end{equation}

Given a prescribed integer $J>0$ and a temporal size $\Delta t>0$, we generate $J$ trajectories, i.e. sampling points, from each ${\bf x}_i\in {\cal S}$ 
by letting 
\begin{equation}\label{generate}
    {\bf x}_{i,j}= {\bf x}_i + \varepsilon_t,  i=1,\ldots,I, j = 1,\ldots,J,
\end{equation}
where $\varepsilon_t$ is some stochastic variable which approximates the random variable $\int_{t}^{t+\Delta t} d X_t$.  For instance, when $X_t=B_t$, a standard Brownian motion, we may choose 
\[
    \varepsilon_t=\sqrt{\Delta t} {\cal N}(0,{\cal I}_d),
\]
where ${\cal N}(0,{\cal I}_d)$ is a $d-$dimensional normal distribution. In the following sections of the paper, $X_t$ might be specified as other stochastic processes which may or may not relate to the PDE \eqref{PDE}. 

The admissible set of sampling points is then defined as 
\[
    \mathcal{S'}=\{{\bf x}_{i,j}|i=1,\ldots,I,j = 1,\ldots,J\}.
\]

\subsection{Error indicators}
To select a set of training points  from an admissible set ${\cal S}'$, 
we need a so-called {\it error indicator} function  which has the following two features. First, it should be more or less equivalent to the error function $e=|u-u_{\theta}|$, namely,  $Ind(x)\sim e({\bf x})$ for all ${\bf x}\in \Omega$. The reason behind this equivalence is as below: when the error $e$ is relatively small (e.g., smaller than a given error tolerance) at a point ${\bf x}$, $u_\theta$  simulates $u$  well at the point ${\bf x}$,  $u_\theta$ has been sufficiently trained at  ${\bf x}$; if $e({\bf x})$ is relatively large, $u_\theta$ does not simulate $u$ well at the point ${\bf x}$, and thus ${\bf x}$ should be chosen as a training point for further training. Secondly, the error indicator should be calculable for each point ${\bf x}$ so that we can use it to select  training points. For this reason,  $e$ itself can not be chosen as an indicator since $u$ is unknown in the process of the training. Consequently, we often need to use the information of the trained $u_\theta$ and  the to-be-solved PDE to construct an appropriate error indicator. For example, $|\mathcal{Q}(u_\theta)|$, the residual of ${\cal Q}$  is a well-known candidate of error indicator which has been used in adaptive sampling method\cite{DeepXDE, EI-PINN, FI-PINN, ADN, DAS-PINN, Selectnet, Rang, ASS-PINN, RAD, ADLGM, RFM}. 
Since  the calculation of the residual  often involves the calculation of some derivatives and  no theory  guarantees the residual equivalents to the error $|u-u_{\theta}|$, in this paper, we will introduce some novel error indicators, in which $u$ is approximated by a so-called {\it empirical value} which is related to $u_\theta$ and the coefficients of the target PDE. The details of our new error indicators  will be presented in Sections 3, 4, and 5.

\subsection{Selection of sampling points}
In this subsection, we explain how to select $I$  to-be-trained sampling points from $\mathcal{S'}$ according to a prescribed indicator  ${\it Ind}$. 

Since the indicator  ${\it Ind}$ is more or less equivalent to the true error $e$,  smaller the value of  ${\it Ind}({\bf x}_{i,j})$ is, better trained  at ${\bf x}_{i,j}$ the NN function $u_\theta$ is. Therefore, we may sort all ${\it Ind}({\bf x}), {\bf x}\in {\cal S}'$ from the largest value to the smallest value, and then select the first $I$ ${\it Ind}$ larger points to generate ${\cal S}''$, the novel set of training points. We call this type of selection the {\it global selection}. Alternatively, we may also generate ${\cal S}''$ by selecting from the $I$ sets of sampling points $\mathcal{S'}_i=\{{\bf x}_{i,j}| j=1,\ldots, J\}$ a point ${\bf x}_i''$ such that 
\[
Ind({\bf x}_i'') =\max \{Ind({\bf x}_{i,j}) : {\bf x}_{i,j}\in\mathcal{S'}_i\}.
\]

Then we let the set ${\cal S}''=\{{\bf x}_i''|i=1,\ldots ,I\}$.  We call this  type of selection as the {\it local selection}.

In summary, to generate ${\cal S}_n$, the set of training points in the $n-$th step, with $n\ge 1$, we first use ${\cal S}_{n-1}$ and a stochastic process to generate the set of admissible training points ${\cal S}'_{n}$. Then by using a properly defined error indicator, we select $I$ training points which have relatively larger indicator values.  
We summarize the procedure of generating ${\cal S}_n$ from ${\cal S}_{n-1}$ by the algorithm described in \ref{adaptive trajectory sampling}.  Note that by randomly sampling ${\cal S}_0$ and iteratively letting ${\cal S}_n= ATS({\cal S}_{n-1})$, we actually generate a sequence of training sets ${\cal S}_n, n=0,1,2,\ldots, N-1$ with $\#{\cal S}_n=I$ for all $n$. Moreover, if in $\Omega$, we connect ${\bf x}_i$ in the $(n-1)$-th step with ${\bf x}_i$ in the $n-$th step by a line segment, we actually obtain $I$ trajectories in $\Omega$.
\SetAlgoCaptionSeparator{}
\label{adaptive trajectory sampling}
\begin{algorithm}
  \renewcommand{\thealgocf}{}
  \caption{$\mathcal{S''}$=ATS($\mathcal{S}, I, J$, Indicator type, Selection type)}
  \SetKwInOut{KwIn}{Initialization}
  \KwIn{$\mathcal{S'}=\emptyset, \mathcal{S'}_i=\emptyset, \mathcal{S}''=\emptyset$.}
  \For{$1\leq i\leq I$}{
        Generate ${\bf x}_{i,j}, j = 1,\ldots,J$ from ${\bf x}_i\in \mathcal{S}$ by \eqref{stochastic process}.\\
        $\mathcal{S'}_i=\mathcal{S'}_i\cup\{{\bf x}_{i,j}\}$.\\
        }
        $\mathcal{S}' = \bigcup\limits_{i=1}^{I}\mathcal{S'}_i$.\\
        Calculate $Ind({\bf x})$ for all ${\bf x}\in \mathcal{S'}$.\\
   
       \For{$1\leq i\leq I$}{
       \If {``Selection type"==``global"}{
            Choose ${\bf x}_i''$ from $\mathcal{S'}$ such that
             \[
             Ind({\bf x}_i'')=\max\left \{ {Ind}({\bf x}): {\bf x}\in \mathcal{S'}\setminus \{{\bf x}''_1,\ldots,{\bf x}''_{i-1} \}\right\}.
             \]
         }
         \If {``Selection type"==``local"}{
            Choose ${\bf x}_i''$ from $\mathcal{S'}_i$ such that
             \[
             Ind({\bf x}_i'')=\max\{ {Ind}({\bf x}_{i,j}): {\bf x}_{i,j}\in\mathcal{S'}_i\}.
             \]
         }
        $\mathcal{S}''=\mathcal{S}''\cup \{{\bf x}''_i\}.$
    } 
\end{algorithm}


\section{The adaptive derivative-free-loss methods(ATS-DFLMs)}\ 

In this section, we apply the ATS technique to the derivative-free-loss method introduced in \cite{Derivative} to construct the adaptive derivative-free-loss methods(ATS-DFLMs).

To present the basic idea of our method, we take ${\cal Q}$ in \eqref{PDE} as the elliptic differential operator defined by 
\begin{equation}\label{ePDE}
    \mathcal{Q}(u):=\frac{1}{2}{\text Tr}(\sigma\sigma^{T} {\text Hess}_x u) + F\cdot\nabla u - G 
\end{equation}
where $F=F(x,u(x)) \in {\mathbb R}^d$, $G=G(x,u(x)) \in  {\mathbb R}^d$ and $\sigma = \sigma (x)\in {\mathbb R}^{d\times d} $ are all known functions, ${\text Hess}_x u$ is the Hessian matrix of $u$ with respect to $x$.

\subsection{The derivative-free loss method(DFLM)}
 The so-called {\it derivative-free loss} in \cite{Derivative} is designed as below. First, we transform the deterministic PDE ${\mathcal Q}(u)=0$ to the Bellman equation 
\begin{equation}\label{bellman}
    u({\bf x})=\mathbb{E}\left[u(X_t) - \int\limits_{0}^{t}G(X_s, u(X_s)) ds\Big|X_0 = {\bf x}\right],
\end{equation}
where $\{X_t\} \in {\mathbb R}^d$ is a $d$-dimensional stochastic process satisfying the  equation
\begin{equation}\label{FSDE}
	dX_t = F(X_t,u(X_t))dt + \sigma(X_t)dB_t,
\end{equation}
for the reason why we can do this transformation, we refer to e.g. \cite{Derivative}.

Secondly, we construct the loss function based on the Bellman equation\eqref{bellman}. We have a set of training points $\mathcal{S}_n=\{{\bf x}_i|i = 1,\ldots, I\}$ at $n$-th step. At the initial step, ${\cal S}_0$ is randomly generated according to a certain distribution. For $n\ge 1$, ${\cal S}_n$ is obtained after the $(n-1)$-th training step.
Now we fix a temporal size $\Delta t>0$, and from ${\bf x}_i$, $i = 1,\ldots, I$, we sample $J_1$($J_1>0$ is a prescribed integer) trajectories, i.e. novel points,  according to the stochastic process $X_t$ determined by \eqref{FSDE}. That is, we let
\begin{equation}\label{disfsde}
{\bf x}_{i,j}= {\bf x}_i + F({\bf x}_i,u_\theta({\bf x}_i))\Delta t + \sigma({\bf x}_i) \sqrt{\Delta t}\mathcal{N}(0,\mathcal{I}_d), i = 1,\ldots,I, j = 1,\ldots,J_1.
\end{equation}

We let the {\it reward} of the point ${\bf x}_i$ be
\begin{equation}
\label{Rij}
    R_{i,j} = G\left({\bf x}_i,u_\theta\left({\bf x}_i\right)\right)\Delta t,
\end{equation}
 and let the {\it empirical value}'s average at the point ${\bf x}_i$ be
\begin{equation}\label{yi}
    y({\bf x}_i) = \frac{1}{J_1}\sum\limits_{j=1}^{J_1}\Big[u_\theta({\bf x}_{i,j})-R_{i,j}\Big].
\end{equation}

The {\it interior} loss function is then defined as 
\begin{equation}\label{DFLMlossin}
    \mathcal{L}^{\Omega}(\theta) := \frac{1}{I}\sum\limits_{i=1}^I\Big(y({\bf x}_i) - u_\theta({\bf x}_i)\Big)^2.
\end{equation}

The interior loss \eqref{DFLMlossin} involves no derivative of the NN function $u_\theta$, and that is why it is called a {\it derivative-free loss}.

{\bf Remark 3.1} The right-hand side of the definition  \eqref{Rij} is not related to the index $j$. Basically, $R_{i,j}$ is an approximation of the integral $\int\limits_{0}^{t}G(X_s, u(X_s)) ds$, where the path generated by the process $X_t$ is the segment starting from the point ${\bf x}_i$ and ending at ${\bf x}_{i,j}$. Therefore, we may also use the definition 
\begin{equation}\label{Rij2}
\begin{aligned}
    R_{i,j} &= \Big[G\left({\bf x}_i,u_\theta({\bf x}_i)\right) + G\left({\bf x}_{i,j},u_\theta({\bf x}_{i,j})\right) \Big]\frac{\Delta t}{2},
\end{aligned}
\end{equation}
to improve the accuracy of the interior loss.

{\bf Remark 3.2}
In the above definition, since ${\bf x}_{i,j}$ is generated by \eqref{FSDE}, then it is possible for ${\bf x}_{i,j}$ to exit the boundary. To overcome this problem, one often needs to use {\it drag back strategy}, e.g., assume that ${\bf x}_i$ are in the $\Omega$ and ${\bf x}_{i,j}$ are points generated from ${\bf x}_i$ using \eqref{FSDE}. If ${\bf x}_{i,j}$ exits $\Omega$, we replace ${\bf x}_{i,j}$ with the intersection of the line $\overline{{\bf x}_i{\bf x}_{i,j}}$ with the boundary $\partial\Omega$, which is denoted as $\hat{{\bf x}}_{i,j}$. Since $\hat{{\bf x}}_{i,j}$ still belong to the $\Omega$ closure, they can participate in the training.

The process $X_t$ depends on the functions $F(\cdot,\cdot)$ and $\sigma(\cdot)$ in the PDE \eqref{ePDE}. So in the case that  the  drift $F$ is nontrivial,  the sampling points ${\bf x}_{i,j}$ might be an imbalance in the neighborhood of ${\bf x}_i$. To avoid the influence caused by $F$, we may prefer to sample points directly using the Brownian motion. For this purpose, we  may use the Cameron-Martin-Girsanov(CMG) theorem\cite{CMG} to transform \eqref{ePDE} to the Bellman equation as below
\begin{equation}\label{Bellman2}
    u({\bf x}) = \mathbb{E}\left[\Big(u(B_{t}) - \mathcal{R}(t)\Big)\cdot\mathcal{D}( t)\Big| B_0={\bf x}\right],
\end{equation}
where 
\begin{equation}\label{u-r}
\mathcal{R}(t) = \int\limits_{0}^{t}\frac{G(B_s,u(B_s))}{\sigma(B_t)\sigma(B_t)^{\top}}ds,\ \ \ \mathcal{D}(t)=exp \left(\int\limits_{0}^{t} \frac{F(B_{s}, u(B_s)}{\sigma(B_t)\sigma(B_t)^{\top}}dB_{s} - \frac{1}{2}\int\limits_{0}^{t}\left\|\frac{F(B_{s}, u(B_{s})}{\sigma\left(B_{s}\right) \sigma\left(B_{s}\right)^{\top}}\right\|^{2}ds\right).
\end{equation}

The interior loss based on \eqref{Bellman2} can be described below. Supposing that we have sampled $I$ points ${\bf x}_i, i = 1,\ldots, I$,  from each point ${\bf x}_i$, we generate $J_1$ novel points by letting
\begin{equation}\label{FSDE2}
    {\bf x}_{i,j}= {\bf x}_i + \sqrt{\Delta t}\mathcal{N}(0,\mathcal{I}_d), i = 1,\ldots ,I, j = 1,\ldots,J_1.
\end{equation}

Moreover, we define
 \begin{equation}\label{CRij}
 \begin{aligned}
     {\mathcal R}_{i,j} &= \frac{G\left({\bf x}_i,u_\theta({\bf x}_i)\right)}{\sigma({\bf x}_{i})\sigma({\bf x}_{i})^{\top}}\Delta t, \ \ \ {\mathcal D}_{i,j} =exp\left(Y_{i,j} \right),
 \end{aligned}
 \end{equation}
where 
 \begin{equation}\label{Yij}
  Y_{i,j} = \frac{F\left({\bf x}_i,u_\theta({\bf x}_i)\right)}{\sigma({\bf x}_{i})\sigma({\bf x}_{i})^{\top}}({\bf x}_{i,j}-{\bf x}_i) -\frac12 \left\|\frac{F\left({\bf x}_i,u_\theta({\bf x}_i)\right)}{\sigma({\bf x}_{i})\sigma({\bf x}_{i})^{\top}}\right\|^2 \Delta t.
\end{equation}

We denote the {\it empirical value}'s average of $u_\theta$ at the point ${\bf x}_i$ by
\begin{equation}\label{yi2}
    y({\bf x}_i) = \frac{1}{J_1}\sum\limits_{j=1}^{J_1}\left[u_\theta({\bf x}_{i,j})-{\mathcal R}_{i,j}\right]{\mathcal D}(i,j).
\end{equation}

Then, we  define the interior loss also by \eqref{DFLMlossin}.

{\bf Remark 3.3} In the calculation of the interior loss on the standard Brownian motion, we may also replace the ${\mathcal R}(i,j)$ 
defined in \eqref{CRij} by 
\begin{equation}\label{CRij2}
\begin{aligned}
    {\mathcal R}_{i,j} &= \left[\frac{G\left({\bf x}_i,u_\theta({\bf x}_i)\right) }{\sigma({\bf x}_i)\sigma({\bf x}_i)^{\top}}+\frac{G\left({\bf x}_{i,j},u_\theta({\bf x}_{i,j})\right)}{\sigma({\bf x}_{i,j})\sigma({\bf x}_{i,j})^{\top}}\right]\frac{\Delta t}{2}, 
\end{aligned}
\end{equation}
and replace the definition \eqref{Yij} by 
\begin{equation}
\begin{aligned}
    Y(i,j) &= \frac{{\bf x}_{i,j}-{\bf x}_i}{2}\cdot\left(\frac{F\left({\bf x}_i,u_\theta({\bf x}_i)\right) }{\sigma({\bf x}_i)\sigma({\bf x}_i)^{\top}}+\frac{F\left({\bf x}_{i,j},u_\theta({\bf x}_{i,j})\right)}{\sigma({\bf x}_{i,j})\sigma({\bf x}_{i,j})^{\top}}\right)\\ 
    &-\left( \left\|\frac{F\left({\bf x}_i,u_\theta({\bf x}_i)\right)}{\sigma({\bf x}_i)\sigma({\bf x}_i)^{\top}}\right\|^2+\left\|\frac{F\left({\bf x}_{i,j},u_\theta({\bf x}_{i,j})\right)}{\sigma({\bf x}_{i,j})\sigma({\bf x}_{i,j})^{\top}}\right\|^2\right)\frac{\Delta t}{4}.
\end{aligned}
\end{equation}

We close the section with a brief explanation of the so-called {\it boundary loss}. Taking the Dirichlet boundary as an example, we may sample $S>0$ points  according to a certain distribution (e.g. the {\it Uniform} distribution, that is, we let $x_s\sim Unif(\partial\Omega), s=1,\ldots, S$, and construct the boundary loss  as 
\begin{equation}\label{DFLMlossb}
    \mathcal{L}^{\partial\Omega}(\theta) = \frac{1}{S} \sum\limits_{s=0}^{S}\left(\mathcal{B}(u_\theta(x_s))-g\right)^2.
\end{equation}
and the total loss function is then defined as
\begin{equation}\label{DFLMloss}
    \mathcal{L}(\theta) := \mathcal{L}^{\Omega}(\theta) + \mathcal{L}^{\partial\Omega}(\theta).
\end{equation}

Finally, we optimize the network parameters $\theta$ at each training iteration by minimizing the loss function defined in \eqref{DFLMloss} using gradient descent, which updates the network parameters $\theta$.

\subsection{Error indicators}
In the following, we present several error indicators for the DFLM solution designed using the so-called {\it empirical value}. That is, the indicators are designed by replacing the exact solution $u$ with its {\it empirical value} in the computation of the error.

The first two  empirical values are calculated  according to the Bellman equation \eqref{bellman}. In the previous subsection, we have generated two sets ${\cal S}_n$ and ${\cal S'}_{n,1}=\{{\bf x}_{i,j}|i=1,\dots,I,j=1,\ldots,J_1\}$ at $n$-th step. Then, given an integer $J_2>0$, we generate an admissible set of training points ${\cal S'}_{n,2}=\{{\bf x}_{i,j}|i=1,\dots,I,j=1,\ldots,J_2\}$. Since the ATS doesn't involve the calculation of the loss, we have the flexibility to choose the values of $J_1$ and $J_2$. To optimize computational efficiency and speed, we choose ${\cal S'}_{n,2}={\cal S'}_{n,1}$ and simplified ${\cal S'}_{n,1}$ to ${\cal S'}_n$. This eliminates the need for additional sampling during the training process and enables us to adaptively select training points from the already obtained set ${\cal S'}_n=\{{\bf x}_i|i=1,\ldots,I*J_1\}$.

Now we just need to construct the error indicators based on the {\it empirical value} of ${\bf x}_i\in{\cal S'}_n$.
Given a temporal size $\Delta t>0$, we generate $J_3$($J_3>0$ is a prescribed integer) sampling points by letting
\[
    {\bf x}_{i,j}= {\bf x}_i + F({\bf x}_i,u_\theta({\bf x}_i))\Delta t + \sigma({\bf x}_i) \sqrt{\Delta t}\mathcal{N}(0,\mathcal{I}_d), i = 1,\ldots,I, j = 1,\ldots ,J_3.
\]
 
We define our first two {\it empirical values} as 
\begin{eqnarray}\label{eav1}
y_1({\bf x}_i) &=& \frac{1}{J_3}\sum\limits_{j=1}^{J_3}\Big[u_\theta({\bf x}_{i,j}) - G({\bf x}_i, u_\theta({\bf x}_i))\Delta t\Big],\\
y_2({\bf x}_i) &=& \frac{1}{J_3}\sum\limits_{j=1}^{J_3}\Big[u_\theta({\bf x}_{i,j}) - \Big[G\left({\bf x}_i,u_\theta({\bf x}_i)\right) + G\left({\bf x}_{i,j},u_\theta({\bf x}_{i,j})\right) \Big]\frac{\Delta t}{2}\Big].    
\end{eqnarray}
    
The third and fourth {\it empirical values} are designed according to the Bellman equation \eqref{Bellman2}. We first generate $J_3$ points by letting 
\[
    {\bf x}_{i,j}= {\bf x}_i + \sqrt{\Delta t}\mathcal{N}(0,\mathcal{I}_d),  j = 1,\ldots J_3.
\]

Let 
 \begin{equation}\label{RjDj1}
 \begin{aligned}
     {\mathcal R}^1_{i,j} &= \frac{G\left({\bf x}_i,u_\theta({\bf x}_i)\right)}{\sigma({\bf x}_i)\sigma({\bf x}_i)^{\top}}\Delta t, \ \ \ {\mathcal D}^1_j =exp\left(\frac{F\left({\bf x}_i,u_\theta({\bf x}_i)\right)}{\sigma({\bf x}_i)\sigma({\bf x}_i)^{\top}}({\bf x}_j-{\bf x}_i) -\frac12 \left\|\frac{F\left({\bf x}_i,u_\theta({\bf x}_i)\right)}{\sigma({\bf x}_i)\sigma({\bf x}_i)^{\top}}\right\|^2 \Delta t \right),
 \end{aligned}
 \end{equation}
and let 
\begin{equation}\label{Rj2}
     {\mathcal R}^2_{i,j}= \left(\frac{G\left({\bf x}_i,u_\theta({\bf x}_i)\right)}{\sigma({\bf x}_i)\sigma({\bf x}_i)^{\top}}+\frac{G\left({\bf x}_{i,j},u_\theta({\bf x}_{i,j})\right)}{\sigma({\bf x}_{i,j})\sigma({\bf x}^{i,j})^{\top}}\right) \frac{\Delta t}{2}, 
 \end{equation}
 and
 \begin{equation}\label{Dj2}
 \begin{aligned}
     {\mathcal D}^2_{i,j} &=exp\Big[\frac{{\bf x}_{i,j}-{\bf x}_i}{2}\cdot\left(\frac{F\left({\bf x}_i,u_\theta({\bf x}_i)\right) }{\sigma({\bf x}_i)\sigma({\bf x}_i)^{\top}}+\frac{F\left({\bf x}_{i,j},u_\theta({\bf x}_{i,j})\right)}{\sigma({\bf x}_{i,j})\sigma({\bf x}_{i,j})^{\top}}\right) \\
    &-\left( \left\|\frac{F\left({\bf x}_i,u_\theta({\bf x}_i)\right)}{\sigma({\bf x}_i)\sigma({\bf x}_i)^{\top}}\right\|^2+\left\|\frac{F\left({\bf x}_{i,j},u_\theta({\bf x}_{i,j})\right)}{\sigma({\bf x}_{i,j})\sigma({\bf x}^{i,j})_{\top}}\right\|^2\right)\frac{\Delta t}{4}\Big].
 \end{aligned}
 \end{equation}
 
We define the third and fourth {\it empirical values} as 
\[
    y_m({\bf x}_i) = \frac{1}{J_3}\sum\limits_{j=1}^{J_3}\Big[u_\theta({\bf x}_{i,j})-{\mathcal R}^{m-2}_{i,j}\Big]{\mathcal D}^{m-2}_{i,j}, m=3,4.
\]

With the above four {\it empirical values}, we define our error indicators at ${\bf x}_i$  by
\begin{equation}\label{Indicator1}
\begin{aligned}
     {\it Ind}_m({\bf x}_i,J_3,\Delta t)= |y_m({\bf x}_i) - u_\theta({\bf x}_i)|, m=1,2,3,4.
\end{aligned}
\end{equation}

Compared to the unknown $u({\bf x}_i)$, the empirical average value $y_m({\bf x}_i)$ is computable during the training process. On the other hand, since the empirical value $y_m({\bf x}_i)$ involves the Bellman equation and thus can be regarded as an approximation of $u({\bf x}_i)$ much better than $u_\theta({\bf x}_i)$. In other words, often we have $|y_m({\bf x}_i) - u({\bf x}_i)|\ll |u({\bf x}_i) - u_\theta({\bf x}_i)|$ and consequently $|y_m({\bf x}_i) - u_\theta({\bf x}_i)|\sim |u({\bf x}_i) - u_\theta({\bf x}_i)|$. This is the reason why we can use $|y_m({\bf x}_i) - u_\theta({\bf x}_i)|$ as our error indicators.

{\bf Remark 3.4} In the simplest case $J_3=1$, to reduce the computational effort, we use ${\bf x}_i\in{\cal S'}_n$ instead of generating a new point ${\bf x}_{i,1}$. 

\subsection{ATS-DFLMs}
With the loss function and error indicators presented previously, we are ready to present our adaptive derivative-free-loss methods(ATS-DFLMs), of which the core is to use the ATS to generate the sets of training points.

Firstly, we define the training steps for the neural network function as $N$, denoted by $n=0,1,\ldots, N-1$.  At the initial step, we randomly generate the set of training points ${\cal S}_0=\{{\bf x}_i|i=1,\dots,I\}$ and then complete the training of $u_\theta$ with ${\cal S}_0$.

Secondly, for the $n$-th step, assuming that we have a set ${\cal S}_n$ obtained from the $(n-1)$-th step.
Then, we use the error indicators $Ind_m$ from the previous subsection to select I to-be-trained points from ${\cal S}_n$. In the ATS-DFLMs for solving \eqref{ePDE}, we can use the {\it global selection} type or {\it local selection} type to generate ${\cal S}_{n+1}$. Since there is no inherent relationship between the training points of two successive steps, we choose the {\it global selection} type to fully explore the domain, such that
\[
    {\it Ind}_m ({\bf x}_{i})=\max \{{\it Ind}_m({\bf x}): {\bf x}\in \mathcal{S'}_n\setminus\{{\bf x}_{1},\ldots,{\bf x}_{i-1}\}\}, 1\le i\le I.
\]
to obtain a new set of training points for the $(n+1)$-th step. 

In summary, for all $n=0,1,\ldots, N-1$, we generate the sets of training points by iteratively using the algorithm
\[
    {\cal S}_{n+1}=ATS({\cal S}_n, I, J_3, Ind_m, ``global") 
\]
and train the NN function $u_\theta$ with it. 

{\bf Remark 3.5} In practice, we prefer not to select an excessively large value for $J_3$, as we only require sorting ${\bf x}_{i,j}$ based on a rough estimate of ${\it Ind}_m({\bf x}_{i,j})$. Through experimental results, we found that the highest accuracy can be achieved and the computation time is reduced when $J_3=1$.

\section{The adaptive PINN method(ATS-PINN)}\

In this section, we apply the ATS technique to the well-known PINN method\cite{PINNs}. 

The PINN can solve many PDEs,  here we only take the elliptic equation \eqref{ePDE} as an example to illustrate our basic idea of the ATS-PINN.  We begin with a brief overview of the PINN (see also \cite{PINNs}). Suppose that we will train $N$ steps the NN function $u_\theta$. For all $0\le n\le N$, the $n$-th step's set of training points in $\Omega$ and on $\partial \Omega$ are denoted by ${\cal S}_n=\{{\bf x}_i|i=1,\ldots,I\}$ and $\mathcal{S}_n^{\partial\Omega}=\{{\bf x}_s|s = 1,\ldots,S\}$, respectively.
The  loss function of the PINN for \eqref{ePDE} is usually defined as :
\begin{eqnarray}
    \mathcal{L}(\theta) =\mathcal{L}^{\Omega}(\theta)+\lambda\mathcal{L}^{\partial\Omega}(\theta),
\end{eqnarray}
where
\begin{eqnarray}
\mathcal{L}^{\Omega}(\theta)=\frac{1}{I}\sum\limits_{i=1}^{I}|\mathcal{Q}(u_\theta({\bf x}_i))|^2,\ \ \mathcal{L}^{\partial\Omega}(\theta)=\frac{1}{S}\sum\limits_{s=1}^{S}|\mathcal{B}(u_\theta({\bf x}_s))-g|^2,
\end{eqnarray}
and $\lambda$ is a weight coefficient dependent or independent of the problem to be solved. 


Next, we explain how to use the ATS to generate the training points in ${\cal S}_n, n=0,\ldots,N$. Usually, ${\cal S}_0$ is randomly sampled from $\Omega$ according to the uniform distribution. Suppose now ${\cal S}_n, n\ge 0$ has been obtained and $u_\theta$ has been trained on ${\cal S}_n$, we generate ${\cal S}_{n+1}$ as below.
For each ${\bf x}_i\in {\cal S}_n,1\le i\le I$, we use a standard Brownian motion to generate $J$ ($J>0$ is a prescribed integer) points. That is, we let  
\[
   {\bf x}_{i,j}= {\bf x}_i + \sqrt{\Delta t}\mathcal{N}(0,\mathcal{I}_d),  j = 1,\ldots,J,
\]
where $\Delta t>0$ is a small prescribed radius. Note that if ${\bf x}_{i,j}$ exits the boundary during the training process, it is replaced with a random sampling point in the domain.

To select $I$ to-be-trained points from $\mathcal{S'}_n=\{{\bf x}_{i,j}|i=1,\ldots,I, j = 1,\ldots,J\}$, we use the following residual-type  error indicator 
\begin{equation}\label{ind-PINN}
     {\it Ind}_P({\bf x}) = |\mathcal{Q}(u_\theta({\bf x}))|, {\bf x}\in \Omega.
\end{equation}

Since the  calculation of  ${\it Ind}_P$ involves to calculate the derivative of $u_\theta$, usually we need to use the {\it automatic differential module} which may lead to a high computational cost.  Therefore, in particular, for high-dimensional PDEs, it is better to use other error indicators to replace ${\it Ind}_P$. For example, we may choose the indicator ${\it Ind}_m$ for $m=3$ or $m=4$, where ${\it Ind}_m$ is the error indicator designed in Section 3. 

Similar to the case in the ATS-DFLMs, we may choose the {\it global} or {\it local} selection type to select a set of training points from $\mathcal{S'}_n$. If we use the {\it global selection} type, the points in ${\cal S}_n$ are selected  such that
\[
    {\it Ind}_P ({\bf x}_{i})=\max \{{\it Ind}_P({\bf x}): {\bf x}\in \mathcal{S'}_n\setminus\{{\bf x}_{1},\ldots,{\bf x}_{i-1}\}\}, 1\le i\le I.
\]

In summary for $n=0,\ldots,N-1$, we iteratively use the algorithm 
\[
    \mathcal{S}_{n+1}=ATS(\mathcal{S}_{n}, I, J, {\it Ind}_{P}, ``global")
\]
to generate all sets of training points. We call the PINN method with the above ATS technique as the ATS-PINN method.

\section{The adaptive temporal-difference method for forward-backward stochastic equations (ATS-FBSTD)}\

In this section, we apply the ATS technique to the temporal difference method for solving forward-backward stochastic equations (FBSTD) proposed in\cite{FBSTD}.

The FBSTD  solves a  parabolic type PDE using the so-called {\it temporal-difference} (TD) learning to train its equivalent stochastic differential equations.  Precisely, to solve the deterministic quasi-linear parabolic PDE
\begin{equation}\label{PDET}
	\begin{cases}
		\frac{\partial u}{\partial t} +\frac12 {\rm Tr}(\sigma\sigma^{\rm T} {\rm Hess}_x u) +F\cdot \nabla u + f=0, \\
		u(T, \cdot) = g,
	\end{cases}
\end{equation}
we first introduce two stochastic processes 
\begin{equation}\label{connection between FBSDE and PDE}
	Y_t = u(t,X_t), \ \ \ \ Z_t = \nabla u(t,X_t),
\end{equation}
so that \eqref{PDET} can be transformed by the It\^o formula\cite{It^o} to the forward-backward stochastic differential 
equations(FBSDEs)  
\begin{equation}\label{FBSDE}
\begin{cases}
    dX_t &= F(t,X_t)dt + \sigma(t,X_t)dB_t, \ X_0 = \xi \in  {\bf R}^d,\\
    dY_t &= - f(t,X_t,Y_t, \sigma^{\rm T}(t,X_t)Z_t)dt + Z^{\rm T}_t \sigma(t, X_t)dB_t, \ Y_T = g(X_T).
\end{cases}
\end{equation}

To simulate the solution $u$ of \eqref{PDET}, we train an NN function $u_\theta$ using \eqref{FBSDE}. To this end, we first discretize \eqref{FBSDE} by the following Euler-Maruyama scheme : let  $[0,T]$ be divided into $N$ subintervals by  points $0={t_0<t_1<\cdots<t_N=T}$, and in each subinterval $[t_n,t_{n+1}], n=0,1,\ldots, N-1$, let
\begin{equation}\label{FBSDE-em}
\begin{cases}
    X_{t_{n+1}} &= X_{t_n} + F(t_n,X_{t_n})\Delta t_n + \sigma(t_n,X_{t_n})\Delta B_{t_n},\\
    Y_{t_{n+1}} &= Y_{t_n}-
    f\left(t_n,X_{t_n},u(t_n,X_{t_n}),\sigma^{\rm T}(t_n,X_{t_n})Z_{t_{n}}\right) \Delta t_n
    + Z^{\rm T}_{t_{n}} \sigma(t_n,X_{t_n})\Delta B_{t_n},
\end{cases}
\end{equation}
where $\Delta t_n=t_{n+1}-t_n$, $\Delta B_{t_n}\sim {\cal N}(\mathbf{0}, \Delta t_n \mathcal{I}_d)$.

Secondly, we generate ${\cal S}_n$, the set of training points at each step $t_n$. As mentioned before, we fix the cardinality of each ${\cal S}_n$ to be $I$. Then we can write $\mathcal{S}_n = \{(t_n,x^{(i)}_{t_n})\in [0, T]\times\Omega |i=1, \ldots, I\}$. 
In \cite{FBSTD}, ${\cal S}_0$ is randomly sampled from $\Omega$  according to a certain distribution such as the {\it Gaussian} or the {\it Uniform} distribution; and once ${\cal S}_n$, $n\ge 0$ has been obtained, the set ${\cal S}_{n+1}$ is generated according to the first formula of \eqref{FBSDE-em} by using the points in ${\cal S}_n$ as starting points. Note that for each $i=1,\ldots, I$, if we connect $x^{(i)}_{t_n}$ and $x^{(i)}_{t_{n+1}}$ with a line segment for all $0\le n\le N-1$, we obtain one trajectory in $\Omega$. Totally we have $I$ such trajectories.

Thirdly, we construct the loss function at the step $t_n$. For each training point $s_{t_n}^{(i)} = (t_n,x^{(i)}_{t_n})$, the value of $u(s_{t_n}^{(i)})$ can be approximated in two ways: the first way is to directly use the value $u_\theta(s_{t_n}^{(i)})$, the second way is to use the second formula of \eqref{FBSDE-em} to obtain a value  
\begin{equation}\label{updata u}
\begin{aligned}
    \overline{u_\theta}(s_{t_n}^{(i)})&=u_\theta(s_{t_{n+1}}^{(i)}) + R(s_{t_n}^{(i)}),
\end{aligned}
\end{equation}
with 
\[
R(s_{t_n}^{(i)})=f\left(s_{t_n}^{(i)},u_\theta(s_{t_n}^{(i)}),\sigma_{t_n}^{\rm T} Z_{t_n}^{(i)}\right) \Delta t_n - Z_{t_n}^{(i){\rm T}}\sigma_{t_n} \Delta B_{t_n}^{(i)}.
\]

The difference between these two approximate values 
\begin{equation}\label{TDerror}
    E_n^{(i)}=E_{n}^{(i)}(s^{(i)}_{t_{n+1}}, \Delta B_{t_n})=\overline{u_\theta}(s_{t_n}^{(i)})-u_\theta(s_{t_{n}}^{(i)})
\end{equation}
is often called as the {\it TD error}  which somehow  indicates how good $u_\theta$ approximates $u$ at the point $s_{t_n}^{(i)}$ (see in \cite{FBSTD}).  
We define the interior loss function and the total loss function at the $n$-th step below
\begin{eqnarray*}
    \mathcal{L}_{n,ite}(\theta) = \frac{1}{I}\sum_{i=1}^{I} |E_{n}^{(i)}|^2, \ \ \ \mathcal{L}_n(\theta)=\mathcal{L}_{n,ite}(\theta)+ \mathcal{L}_{T}(\theta),
\end{eqnarray*}
where $\mathcal{L}_{T}(\theta)$ is a loss  associated with the terminal condition (see the second equation of \eqref{PDET} or \cite{FBSTD}).

Note that the above FBSTD method is different from  the so-called {\it Monte Carlo} method\cite{RL} which updates the $u_\theta$ only when all trajectories have been completed. Here for the FBSTD, we construct a loss for each time step $t_n$ so we can improve $u_\theta$ at each step $t_n$, without waiting for the final outcome of each whole trajectory. 

In the following, we explain how to use the ATS technique to generate the sets of training points in the FBSTD. 
Recall that in \cite{FBSTD}, ${\cal S}_{n+1}$ are randomly generated by using the stochastic process  \eqref{FBSDE-em} and the starting points in ${\cal S}_n$. By this kind of random sampling, some selected training points in  ${\cal S}_{n+1}$ may not need further training since $u_\theta$ has already approximated $u$ very well at those points. To improve the efficiency of training, we should try to select the training points  where $u_\theta$  approximates $u$ not so well. Based on this idea, we design our  ATS for the FBSTD as below.

First, from a starting point $s_{t_n}^{(i)}\in {\cal S}_n$, instead of only generating one point using the first formula of \eqref{FBSDE-em},  we generate $J$ ($J>0$ is a prescribed integer) admissible points $x_{t_n}^{(i,j)}, j=1,\ldots, J $. We denote  $s_{t_{n+1}}^{(i,j)}=(t_{n+1}, x_{t_n}^{(i,j)})$ and let $\mathcal{S}^{(i)}_n = \{s_{t_{n+1}}^{(i,j)} |j=1,\ldots,J\}$.
Then we explain how to select one to-be trained point from $\mathcal{S}^{(i)}_n$. To this end, we choose the TD error defined in  \eqref{TDerror} as our error indicator. That is, we let  
\begin{equation}\label{ind-FBSTD}
 {\it Ind}_F(s_{t_{n+1}}^{(i,j)})=\left|E_{n}^{(i)}(s_{t_{n+1}}^{(i,j)}, \Delta B_{t_n}^{(i,j)})\right|, \ \ \ j=1,\ldots, J. 
\end{equation}

With this indicator, we select a point $s_{t_{n+1}}^{(i)}$ from  $\mathcal{S}^{(i)}_n$ such that
\[
    {\it Ind}_F (s_{t_{n+1}}^{(i)})=\max \{{\it Ind}_F (s): s\in \mathcal{S}^{(i)}_n\}.
\]

Repeat the above procedure for all $1\le i\le I$, we obtain all to-be-trained points in the $(n+1)$-th step.  In summary, instead of using \eqref{FBSDE-em} to randomly generate the set ${\cal S}_{n+1}$, here we use the algorithm 
\[
    {\cal S}_{n+1}={\rm ATS}({\cal S}_{n}, I,J, {\it Ind}_F, ``local" )
\]
to generate the $(n+1)$-th step's training points. Here we emphasize that for the ATS-FBSTD, we  use the so-called {\it local selection} type to select the training points in order to maintain the temporal association between two consecutive steps. This is different from the case in the ATS-DFLMs where both the {\it local selection} type and {\it global selection} type can be used. We call the FBSTD method with the above adaptive selection of training points the ATS-FBSTD method.

\section{Numerical Results}\
\label{numerical-experiments}

In this section, we test four numerical examples: the 2D Laplace  equation with a singular solution, the high-dimensional Poisson equation, the Black-Scholes 
equation, and the quadratically growing equation. In all our tests, we use a ResNet with Adam\cite{Adam} optimizer to train the 
approximate solution $u_\theta$.  The accuracy of $u_\theta$  is indicated by the {\it relative error} defined by ${\rm RE}=\frac{\|u_\theta-u\|_{L^2}}{\|u\|_{L^2}}$. 

In the first two examples, we test the performance of the ATS-DFLMs and the ATS-PINN method on a ResNet with 3 residual blocks and 60 neurons per layer; we also set the learning rate $\alpha=0.01$ and the activation function to be SWISH\cite{swish}. Note that we use the $Ind_m, m=2,4$ error indicator whenever possible because the trapezoidal format gives a more accurate value of the empirical values $y_m$.

In the third and fourth examples, we will test the performance of the ATS-FBSTD on a ResNet with 2 residual blocks and 256 neurons per layer; we  choose the learning rate $\alpha=0.005$ and the activation function to be $\sin$.


{\bf Example 1}
We consider the  Laplace equation 
\begin{equation}
\label{singularity}
\begin{cases}
    \Delta u = 0,\ &{\rm in}\ \Omega = \{(r,\theta):0\leq r\leq 1,0\leq \theta \leq \frac{\pi}{6}\},\\
    u = h,\ &{\rm on}\ \partial\Omega,
\end{cases}
\end{equation}
where the boundary function $h$ is chosen such that \eqref{singularity} admits the exact solution
\[
    u(r,\theta) = r^{\frac{2}{3}}\sin\left(\frac{2}{3}\theta\right).
\]
Note that $u$ has a singularity at the origin, see also in \cite{Derivative}.

With this example, we  test the performance of both the ATS-DFLM and the ATS-PINN.  For comparison, we will also solve \eqref{singularity} by the original DFLM and PINN without using the ATS. In each training step of each method, we always choose $I=500$, $S=300$. That is, in each training step, we sample 500 training points in $\Omega$ and 300 training points on $\partial\Omega$. Moreover, we let  $\Delta t=5e-04$ and choose $J_1=20$, $Ind_2$ for ATS-DFLM. For the ATS-PINN, we set $J=10$, $\lambda=1$, and choose the error indicator to be $Ind_P$.  

Presented in Table \ref{tab:numresluts of singularity} are the relative errors and computation times after 20000 training steps of different methods. We first observe that these methods take almost the same time to train 20000 steps, which means the additional ATS subroutine does not take much more computational cost compared to the original deep PDE solvers. Moreover, we observe that compared to the PINN and DFLM, the ATS-PINN and ATS-DFLM remarkably improve the accuracy. Note that for the ATS-DFLM, we test two cases: the case $J_3=2$ and the case $J_3=1$. We observe that the case $J_3=1$ not only leads to high accuracy but also takes a low computation time. Therefore in the rest numerical tests of this section, we always choose $J_3=1$ for the ATS-DFLM. 
\begin{table}[H]
	\centering
      \caption{Relative errors and computation time for Example 1.}
	\renewcommand{\arraystretch}{1.0}
		\begin{tabular}{lcc}
			\hline
			Method       & RE        & Time(s)\\
			\hline
			PINN         & 1.833e-02 & 449\\
			ATS-PINN     & 9.765e-03 & 535\\
			DFLM         & 2.874e-02 & 422  \\
			ATS-DFLM($J_3=2$)     & 1.442e-02 & 641\\
            ATS-DFLM($J_3=1$)    & 6.431e-03& 455\\
			\hline
		\end{tabular}
	\label{tab:numresluts of singularity}
\end{table}

\begin{figure}[H]
\centering
\includegraphics[height=5.6cm, width=8cm]{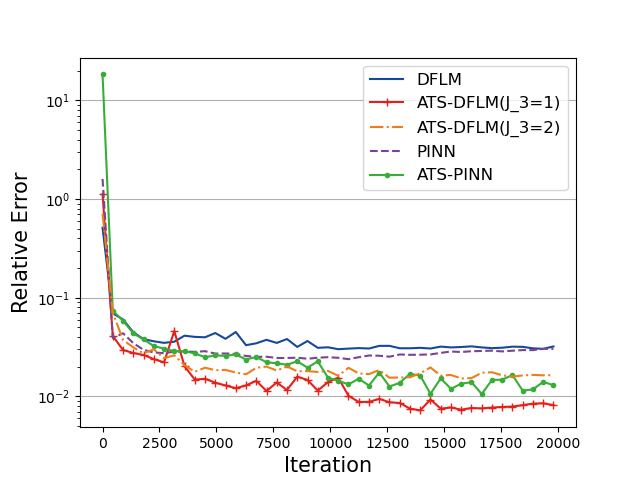}
\caption{Relative errors w.r.t. iteration steps for \eqref{singularity}.}
\label{fig:singularity RE}
\end{figure}

Depicted in Fig \ref{fig:singularity RE} is the dynamic change of the relative errors with respect to the training steps. Again, we observe that the method with the ATS technique performs better than the corresponding method without using the ATS. 

\begin{figure}
\centering
		\subfigure[$0$-th step.]{
        \centering
        \begin{minipage}[t]{0.23\linewidth}
		\centerline{\includegraphics[width=\linewidth]{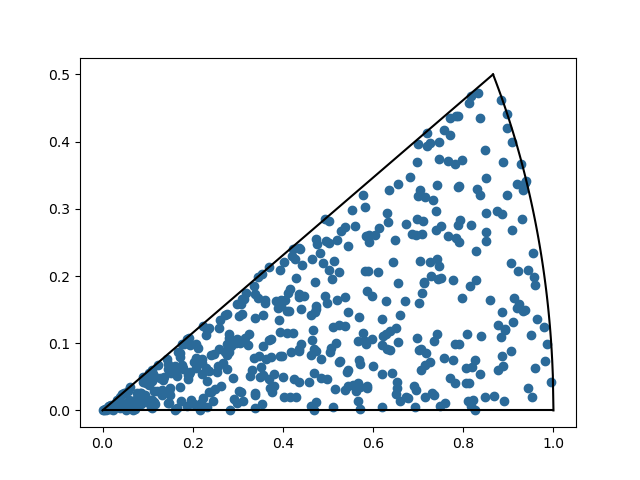}}
	\end{minipage}
        } 
		\subfigure[$2000$-th step.]{
        \centering
        \begin{minipage}[t]{0.23\linewidth}
		\centerline{\includegraphics[width=\linewidth]{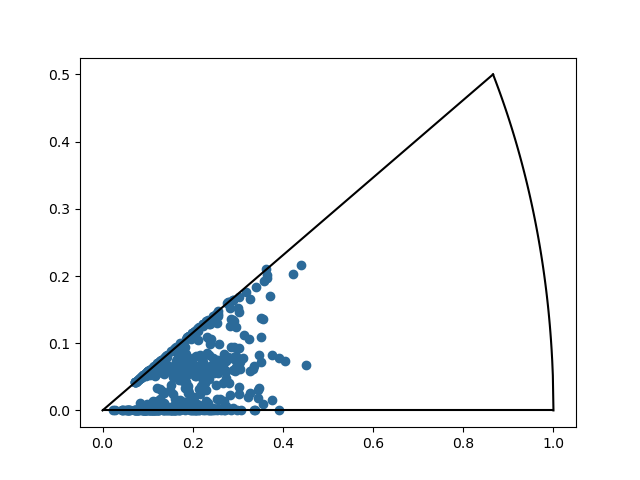}}
	\end{minipage}
        }
        \subfigure[$4000$-th step.]{
        \centering
        \begin{minipage}[t]{0.23\linewidth}
		\centerline{\includegraphics[width=\linewidth]{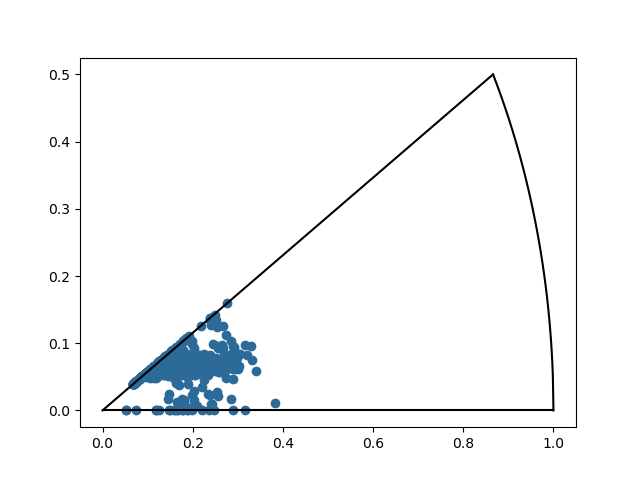}}
	\end{minipage}
        }
        \subfigure[$6000$-th step.]{
        \centering
        \begin{minipage}[t]{0.23\linewidth}
		\centerline{\includegraphics[width=\linewidth]{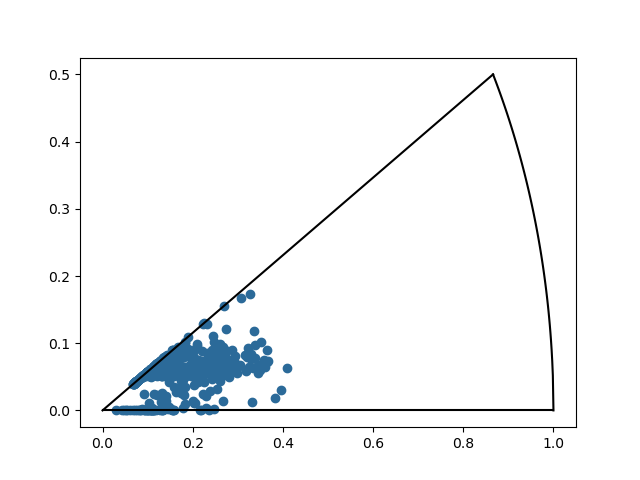}}
	\end{minipage}
        }
 	\caption{The distribution of training points at different steps.}
	\label{fig:transition}
\end{figure}

To show the influence of the ATS on each step's training points, we present in Fig \ref{fig:transition} all 500 interior training points  at $0/2000/4000/6000$-th steps using ATS-DFLM($J_3=1$).  We observe that the training points are uniformly distributed in the initial step, and as the number of steps increases, they gradually move toward the singularity region.

{\bf Example 2} 
We consider the high-dimensional linear Poisson equation 
\begin{equation}
\label{linearpoisson}
\begin{cases}
-\Delta u(x) - G(x) = 0,\ &{\rm in}\ \Omega=[-1, 1]^d,\\
u(x)= h(x), \ &{\rm on}\ \partial\Omega,
\end{cases}
\end{equation}
where $G(x) = \frac1d\left(\sin(\frac1d\sum\limits_{i=1}^{d}x_i)-2\right)$ and the boundary $h(x) = \left(\frac1d\sum\limits_{i=1}^{d}x_i\right)^2 + \sin\left(\frac1d\sum\limits_{i=1}^{d}x_i\right)$.  This problem  admits the exact solution
\[
	u(x)=\left(\frac1d\sum\limits_{i=1}^{d}x_i\right)^2 + \sin\left(\frac1d\sum\limits_{i=1}^{d}x_i\right).
\]

We will test the performance of the ATS-DFLM, ATS-PINN, DFLM, and PINN on the problem \eqref{linearpoisson} with $d=10, 20, 30, 50$. All parameters used in the algorithms for this example are chosen as the same as those for \eqref{singularity}. In particular, we use both $Ind_P$ and $Ind_4$ in the ATS-PINN. When the ATS-PINN method uses the $Ind_4$ error indicator, we still choose $J_3 = 1$.

Table \ref{table:lp re} shows the relative errors and computation times of different methods in different dimensions after 20000 training steps. From this table, we  observe that for all dimensions, the computation time of the DFLM is less than that of the PINN, and the larger the dimension, the greater the  difference. We also observe that for all dimensions,  the ATS-PINN achieves higher accuracy than the PINN does, and the ATS-DFLM achieves higher accuracy than the DFLM does, which implies that our ATS improves the accuracy of the original method independently of the dimension.
We also would like to mention that for the ATS-PINN, the error indicator $Ind_4$ improves the accuracy of the PINN by two orders of magnitude. At the same time, computation time is reduced by half compared to residual type error indicators $Ind_P$.
\begin{table}
	\centering
     \caption{Relative errors and computation times for Example 2.}
	\setlength{\tabcolsep}{0.5mm}{
		\begin{tabular}{cc|cc|cc|cc|cc}
			\hline
			\multicolumn{2}{c|}{\multirow{2}{*}{\diagbox{Method}{Dimension}}}
			&\multicolumn{2}{c|}{d=10} & \multicolumn{2}{c|}{d=20}
            & \multicolumn{2}{c|}{d=30}
            & \multicolumn{2}{c}{d=50}\\
			\cline{3-10}
			\multicolumn{2}{c|}{}&RE&Time(s)&RE&Time(s)&RE&Time(s)&RE&Time(s)\\
			\cline{1-10}
			\multicolumn{2}{c|}{PINN} &5.657e-02&676& 5.756e-02&1223&8.556e-02&1767& 1.269e-01& 2529 \\
   
			\multicolumn{2}{c|}{ATS-PINN($Ind_P$)}&1.376e-02&1638& 2.515e-03&2431&1.219e-02&4163&2.619e-02&6315\\
   
			\multicolumn{2}{c|}{ATS-PINN($Ind_4$)}&7.176e-04&832& 8.354e-04&1324&1.187e-03&2041&1.144e-03&2985  \\
   
            \multicolumn{2}{c|}{DFLM}&9.138e-04&466 & 1.098e-03&643& 3.452e-03&814& 2.430e-03&1324  \\

   		\multicolumn{2}{c|}{ATS-DFLM}&5.982e-04&485&8.535e-04&689&6.854e-04 &818& 1.380e-03&1370 \\
			\hline
		\end{tabular}
	}
	\label{table:lp re}
\end{table}

\begin{figure}[H]
\centering
\subfigure[d=10]{
\begin{minipage}[t]{0.23\linewidth}
\centering
\includegraphics[width=\linewidth]{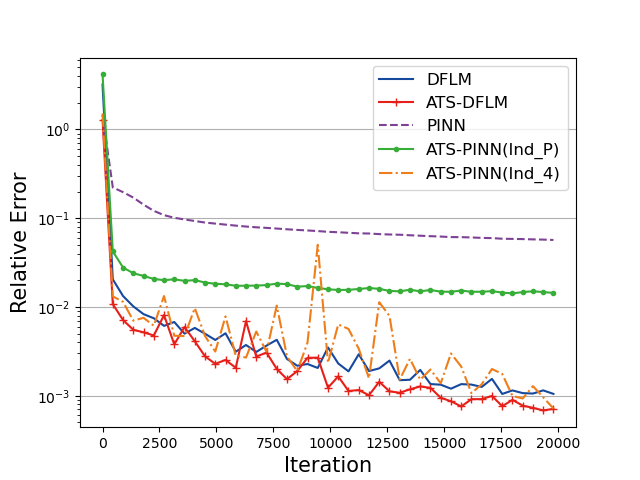}
\end{minipage}%
}%
\subfigure[d=20]{
\begin{minipage}[t]{0.23\linewidth}
\centering
\includegraphics[width=\linewidth]{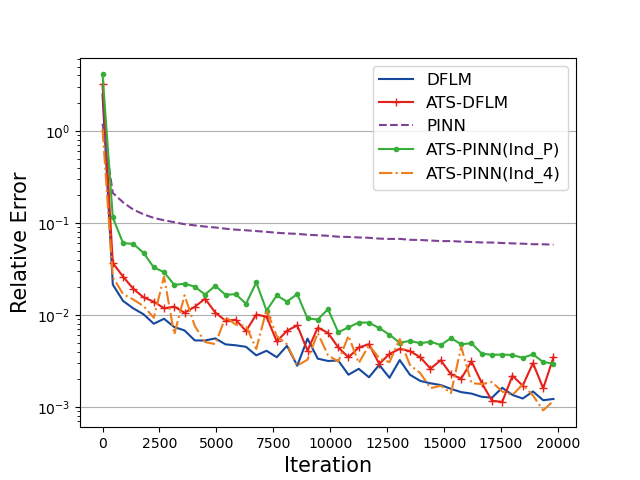}
\end{minipage}%
}%
\subfigure[d=30]{
\begin{minipage}[t]{0.23\linewidth}
\centering
\includegraphics[width=\linewidth]{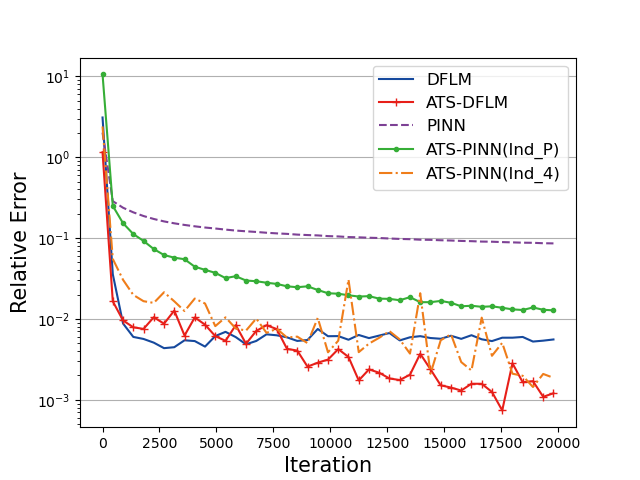}
\end{minipage}
}%
\subfigure[d=50]{
\begin{minipage}[t]{0.23\linewidth}
\centering
\includegraphics[width=\linewidth]{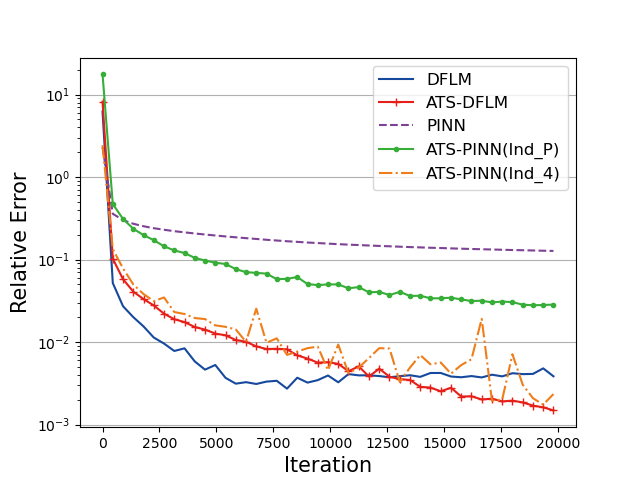}
\end{minipage}
}%
\centering
\caption{The dynamic change of relative errors of \eqref{linearpoisson}.} 
\label{fig:lperror}
\end{figure}

Fig \ref{fig:lperror} shows the dynamic change of the relative errors with respect to the training steps. We observe that the change curves of ATS-DFLM and DFLM are very similar at 10 and 20 dimensions. At 30 and 50 dimensions, although ATS-DFLM does not converge as fast as DFLM until 10000 steps, ATS-DFLM continues to converge after 10000 steps when DFLM has already converged. For different dimensions, we also find that for sufficiently large iteration steps, the ATS-PINN converges to a lower relative error than PINN upfront, where ATS-PINN($Ind_4$) reaches the lowest relative error among PINN and ATS-PINN($Ind_P$), but its curves are more oscillatory.
In this sense, we may say that the performance of the ATS technique is independent of the dimension.

{\bf Example 3}
We consider the following Black-Scholes equation 
\begin{equation}\label{BS}
\begin{cases}
\dfrac{\partial u}{\partial t}(t, x)=-\dfrac{1}{2}{\rm Tr}[0.01\diag(x^2) {\rm Hess}_x u(t, x)]+0.1(u(t, x)-(\nabla u(t, x),x)),\\
u(T,x)= \|x\|^2,
\end{cases}
\end{equation}
which admits the exact solution  $u(x,t)=\exp((0.1+0.1^2)(T-t))\|x\|^2, x\in \mathbb{R}^{100}$. 

We use this example to test the performance of the ATS-FBSTD. For comparison, we  also solve \eqref{BS} by the FBSTD. Recall that in \cite{FBSTD}, we have used two approaches to calculate the derivative of the NN function $u_\theta$: the one uses the {\it automatic differential module} is named as the FBSTD1 and the one uses another NN function to simulate $\nabla u_\theta$ is named as the FBSTD2. In this paper, we will follow these notations and name the FBSTD1 and FBSTD2 with the ATS as the ATS-FBSTD1 and the ATS-FBSTD2, respectively. 

To simulate the solution of \eqref{BS}, we partition  $[0,T]$ into $N=50$ equally sized subintervals and set the starting point ${\bf x}_0 = (1, 0.5,...,1,0.5) \in \mathbb{R}^{100}$. We fix the number of training points to be $I=512$ for all four methods(FBSTD1, FBSTD2, ATS-FBSTD1, ATS-FBSTD2). In addition, for both the ATS-FBSTD1 and the ATS-FBSTD1, we set  $J=10$. For the starting point ${\bf x}_0$, we denote the relative error ${\rm RE}_0=\frac{|u_\theta(0, {\bf x}_0)-u(0, {\bf x}_0)|}{|u(0, {\bf x}_0)|}$. 

Table \ref{tab:numresluts of BS} shows the relative errors RE and ${\rm RE}_0$ of the above methods after 20000 training steps. From the table, we find that the ATS-FBSTD obtains lower relative errors than the FBSTD. Among them, the ATS-FBSTD1 outperforms the others, improving the relative error by an order of magnitude compared to the FBSTD1. We can say that the addition of the ATS technique improves the computational accuracy for solving high-dimensional parabolic equations.  
\begin{table}[H]
	\centering
      \caption{Relative errors and computation time for Example 3.}
	\renewcommand{\arraystretch}{1}
	\setlength{\tabcolsep}{3mm}{
		\begin{tabular}{lccccc}
			\hline
			Method    & $RE$ & $RE_0$ &Time(s) \\
			\hline            
			FBSTD1    & 3.110e-03 & 3.234e-03 & 530 \\
            FBSTD2    & 3.640e-03 & 1.934e-02 & 448 \\
            ATS-FBSTD1    & 7.078e-04 & 6.561e-07 & 546 \\
            ATS-FBSTD2    & 1.500e-03 & 7.185e-05 & 507 \\
            \hline
		\end{tabular}
	}
	\label{tab:numresluts of BS}
\end{table}
\begin{figure}[H]
    \centering
    \includegraphics[width=0.45\linewidth]{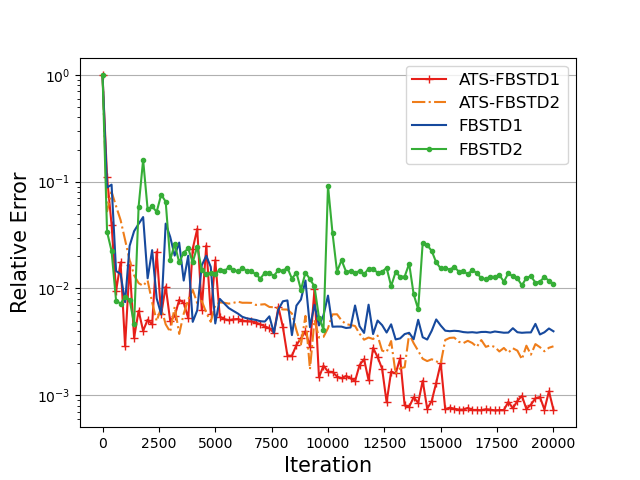}
    \caption{The dynamic change of relative errors of \eqref{BS}.}
    \label{fig:BS-L2}
\end{figure}
Depicted in Fig \ref{fig:BS-L2} is the dynamic change of the relative errors RE with respect to the training steps.  Again, we observe that both ATS-FBSTD1 and ATS-FBSTD2  converge to lower relative errors than  FBSTD1 and FBSTD2. 

\begin{figure}
\centering
    \subfigure[The fit of FBSTD1.]{
        \begin{minipage}[t]{0.45\linewidth}
        \centering
        \centerline{\includegraphics[width=\linewidth]{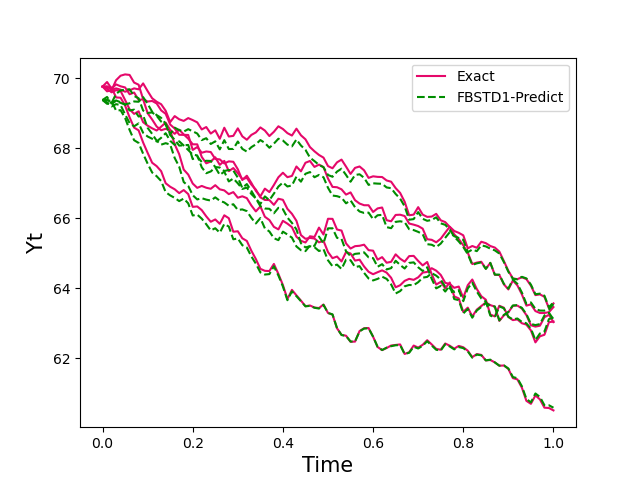}}
		\label{fig:BS-FBSTD-tra}
		\end{minipage}}
    \subfigure[The fit of ATS-FBSTD1.]{
        \begin{minipage}[t]{0.45\linewidth}
        \centering
        \centerline{\includegraphics[width=\linewidth]{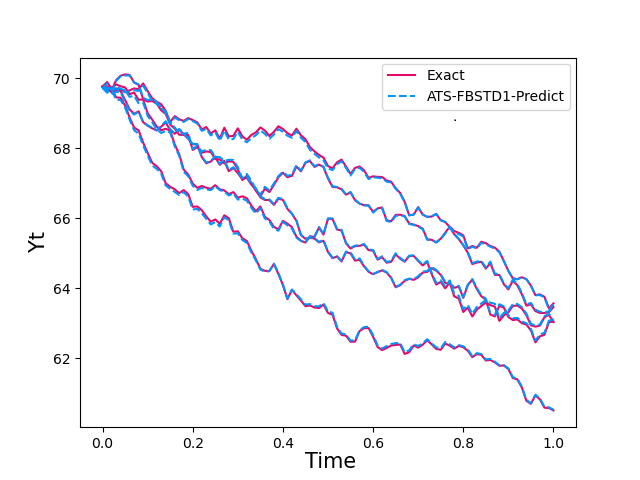}}
		\label{fig:BS-AFBSTD1-tra}
		\end{minipage}}
\caption{The fit of FBSTD1 and ATS-FBSTD1 solution to the exact solution of \eqref{BS}. The red lines represent the exact solution, whereas the green and blue lines represent the FBSTD1 solution and the ATS-FBSTD1 solution respectively.}
\label{fig:BS-tra}
\end{figure}

Fig \ref{fig:BS-tra} shows how  the approximate solutions computed with FBSTD1 and ATS-FBSTD1 fit  the exact solution on 5 randomly generated trajectories. We found that the FBSTD1 solution does not fit the exact solution very well in the beginning, while the ATS-FBSTD1 solution fit the exact solution very well over the whole trajectory, which implies that the ATS-FBSTD1 may have better generalization ability than the FBSTD1.

{\bf Example 4} 
We consider the nonlinear quadratically growing equation
\begin{equation}\label{QG}
	\begin{cases}
		\frac{\partial u}{\partial t}(t, x)+\left \| \nabla u(t,x) \right \|^2 +\dfrac{1}{2}\Delta u(t,x)=h(t,x, u, \nabla u),\\
		u(T,x)=\sin\left(\left[\frac{1}{d}\|x \|^2\right]^\alpha\right),
	\end{cases}
\end{equation}
where $\kappa = 0.2$ and $h(t,x, u, \nabla u)= 2\|x\|^2d^{-2} \left\{\kappa^2 l^{2\kappa-2} [2 \cos^2 l^{\kappa}-\sin l^\kappa] +\kappa(\kappa-1)l^{\kappa-2} \cos(l^{\kappa})\right\}$, with $l=l(t, x) = T-t+\frac{1}{d}\|x\|^2$. This equation admits the analysis solution
\[
	u(t,x)=\sin\left(\left[T-t+\frac{1}{d}\| x \|^{2}\right]^{\kappa}\right).
\]

We will solve \eqref{QG} also by the FBSTD and the ATS-FBSTD. To this end,  we partition  $[0,T]$ into $N=30$ equally sized subintervals and set the starting point ${\bf x}_0 = (0, 0,...,0,0) \in \mathbb{R}^{100}$. We fix the number of training points to be $I=512$ for all four methods and $J=10$ for the ATS-FBSTD. Since the equation \eqref{QG} is nonlinear, the coefficient $F$ in its corresponding FBSDEs \eqref{FBSDE} depends on $u$. Therefore when we generate ${\cal S}_{n+1}$ from ${\cal S}_n$  by using the first equation of \eqref{FBSDE}, we should replace $u$ in $F(\cdot, u)$ by $u_\theta$ which has been trained in the previous step.

Table \ref{table:NumResluts of QG} presents the detailed relative errors RE and ${\rm RE}_0$ of different methods after 20000 iteration steps. We find  that for this nonlinear problem, the ATS-FBSTD achieves better computational accuracy than the FBSTD does.
\begin{table}[H]
	\centering
	\caption{Relative errors and computation time for Example 4.}
		\renewcommand{\arraystretch}{1}
			\setlength{\tabcolsep}{3mm}{
				\begin{tabular}{lccccc}
					\hline
					Method    & $RE$   &$RE_0$  & Time(s) \\
					\hline
					FBSTD1   & 3.845e-03 & 6.159e-05 & 654  \\
					FBSTD2    & 5.500e-03 & 5.611e-03 & 558   \\
                    ATS-FBSTD1   & 2.513e-03 & 2.054e-06 & 692  \\
					ATS-FBSTD2   & 3.385e-03& 2.327e-05 & 589  \\
					\hline
				\end{tabular}
			}
	\label{table:NumResluts of QG}
\end{table}



\bigskip
\bigskip

{\bf Credit authorship contribution statement}
{Xingyu Chen:} Methodology, Software, Writing
 the original draft
 {Jianhuan Cen:} Software 
 {Qingsong Zou:} Conceptualization, Methodology, Writing the original draft 

{\bf Declaration of competing interest}
The authors declare that they have no known competing ﬁnancial interests or personal relationships that could have appeared to inﬂuence the work reported in this paper. 

{\bf Acknowledgements} 
The research was supported in part by the National Key R\&D Program of China (2022ZD0117805), the NSFC Grant 12071496, the Guangdong Provincial NSF 2023A1515012097.



\end{document}